\documentclass[12pt,a4paper]{article}
\usepackage[utf8]{inputenc}
\usepackage{amsmath}
\usepackage{amsfonts}
\usepackage{amssymb}
\usepackage{makeidx}
\usepackage{graphicx}
\usepackage{moreverb}

\newtheorem{theorem}{Theorem}

\usepackage{float}
\begin{document}
	\title{Improved Estimator of the Conditional Tail Expectation 
		in the case of heavy-tailed losses}
	
	\author{Mohamed Laidi\\National High School of Technology, Algiers, Algeria.\\ LRDSI laboratory, Blida 1 university, Blida, Algeria.\\Email: mohamed.laidi@enst.dz
	\and 
	Abdelaziz Rassoul\\High School of Hydraulics, Blida, Algeria.\\Email: a.rassoul@ensh.dz
	\and  
	Hamid Ould Rouis \\ University of Blida 1, Blida, Algeria.\\(Email: houldrouis@hotmail.com)
	}
\date{}
\maketitle
\begin{abstract}
In this paper, we investigate the extreme-value methodology, to propose an improved estimator of the conditional tail expectation ($CTE$) for a loss distribution with a finite mean but infinite variance.
The present work introduces a new estimator of the $CTE$ based on the bias-reduced estimators of high quantile for heavy-tailed distributions. The asymptotic normality of the proposed estimator is established and checked, in a
simulation study. Moreover, we compare, in terms of bias and mean squared
error, our estimator with the known old estimator.\\
\end{abstract}



\section{Introduction and Motivation}\label{sect-1}
Risk management is a subject of concern in finance and actuarial science. Protecting against financial and actuarial risks is essential in order to anticipate financial crises or major insurance claims. For this, there are many tools to quantify and predict risk (risk measures). They make it possible to evaluate a level of danger of a risk, but also to compare different risks between them. Quantification, prevention, comparison and prediction of risk are essential elements of our society.\\     
One of the best known and used risk measure is the Value-at-Risk (or $VaR$), it is introduced in the 1990s by Morgan \cite{morgan1997creditmetrics} (see Jorion \cite{jorion1997value} for more details). The idea of the $VaR$ is as follows: we fixe a threshold $ \alpha $ and calculate a value $VaR(\alpha)$ which will be such that the probability that the catastrophe will occur is smaller than $ \alpha $. $VaR(\alpha)$ can be seen as the amount of extra capital and business needs in order to reduce the probability of going bankrupt to $\alpha$.\\
In statistical terms, the $VaR$ of level $\alpha \in \left] 0,1\right[ $ corresponds to the level $(1-\alpha)$ quantile of the distribution function of the losses. This risk measure has several flaws. It gives only one point information to the quantile $q(\alpha)$ and no information beyond this point; it does not take into account the importance of the disaster when it occurs but only its frequency.\\
To respond to the need for theoretical practical principles,
Artzner et al. 1999 \cite{artzner} introduced the concept of coherent risk measure. A risk measure is said to be coherent if it satisfies the four properties: monotony, positive homogeneity, invariance by translation and subadditivity, (for more details see, \cite{artzner}).\\       
Unfortunately, the $VaR$ risk measure is not subadditive which implies that it is not coherent.
For this reason, many authors proposed to replace standard $VaR$ with alternative risk measures such as Conditional Tail Expectation ($CTE$), also called in the literature on risk measure  Tail Value-at-Risk ($TVaR$) or Expected Shortfall ($ES$) (see, for example, \cite{kaas2008modern}, \cite{dhaene2008}, \cite{acerbi2004coherent}, \cite{acerbi2002expected}, \cite{acerbi2002coherence}, \cite{rockafellar2000optimization}, \cite{rockafellar2002conditional}, \cite{wirch1999synthesis} and \cite{yamai2002comparative}). 

This risk measures gives information of the distribution beyond the $VaR(\alpha)$. Thus, unlike the $VaR$, it takes into account the information contained in the tail of the distributi which is essential in the finanacial context. The $CTE$ has been studied by many authors such that: \cite{artzner}, \cite{cai2007optimal}, \cite{cai2015estimation}, \cite{tasche2002expected}, \cite{brazauskas2008estimating}, \cite{landsman2003}.\\
 The application of these alternative measures has gained interest growing in literature and industry. In particular, the $CTE$, due to its properties and its effectiveness in different fields such as finance and actuaries, its use and studies keep increasing (see, e.g., \cite{dhaene2008}, \cite{tan2011}, \cite{hakim2018pot}, \cite{el2019optimality} and references therein).\\
In the literature, the estimation of the $CTE$ and its asymptotic normality under the assumption that the second moment of the loss variable is finite has been established but this assumption is very restrictive in practical problems. For this reason, Necir et al (2010) \cite{necir2010} has extended that in the case of infinite variance. In this paper, we propose an improvement of the estimator established by Necir et al. (2010) \cite{necir2010}, our considerations are based on the bias-reduced estimators of high quantile for heavy-tailed distributions introduced by Li (2010) \cite{li2010}, also we show its efficiency and its asymptotic normality theoretically, finally we prove the performance of our estimator by some results of simulation study. For this aim, we proceed as follow :\\
 Let $X$ be a loss random variable with cumulative distribution function
(cdf) $F$. $F$ is assumed to be continuous throughout the present paper and  defined on the entire real line, with negative loss interpreted as gain. The $CTE$ of the risk $X$ is then defined, for every $t\in (0,1)$, by the function $ \mathbb{C}(t) $ writed as follow
\begin{equation}
\mathbb{C}(t)=\mathbf{E}\left( X|\,X>\mathbf{Q}(t)\right) ,  \label{eq 4}
\end{equation}%
where $\mathbf{Q}(t)=\inf \left\{ x:F(x)\geq t\right\} $ is the quantile
function corresponding to the cdf $F$. Since $F$ is continuous, we
easily check that $\mathbb{C}(t)$ is equal to 
\begin{equation}
\mathbb{C}(t)=\frac{1}{1-t}\int_{t}^{1}\mathbf{Q}(s)ds.  \label{eq 5}
\end{equation}%

In this paper we are interesting by the heavy-tailed distribution. A model $F $ is said to be heavy-tailed, with a tail index $\alpha $ $\left( \alpha >0\right) $, if the tail function $\overline{F}=1-F\in \mathcal{RV}_{-\alpha
}$, where denotes the class of regularly varying functions with index of regular variation equal to $%
\alpha $, i.e., non-negative measurable functions $g$ such that, for all $%
x>0 $, $g(tx)/g(t)\rightarrow x^{\alpha }$, as $t\rightarrow \infty $
(see, \cite{gnedenko1943} for more details).

Suppose that $X_{1},X_{2},...,X_{n}$ are independent and identically distributed ( i.i.d.) random variables with common
distribution function $F$ which has regularly varying tails with index $%
\alpha >1$, i.e.,
\begin{equation}
\underset{t\rightarrow \infty }{\lim }\frac{1-F(tx)}{1-F(t)}=x^{-\alpha },%
\text{ for every }x>0.  \label{eq1-1}
\end{equation}%
This class includes a number of popular distributions such as Pareto,
generalized Pareto, Burr, Fr\'{e}chet, Student, etc., which are known to be
appropriate models for fitting large insurance claims, fluctuations of
prices, log-returns, etc. (see, e.g., \cite{beirlant2001}). In this paper, we restrict ourselves to this class of distributions. For more
information on the topic and, generally, on extreme value models and their
manifold applications, we refer to the monographs by  \cite{deHaan2006}, \cite{peng1998asymptotically} and \cite{resnick2007}.\\
In particular, if the index $\alpha \in \left( 1,2\right) $, then, condition
(\ref{eq1-1}) is equivalent to the statement that $F$ has infinite second
moment. In this case, in order to estimate the $CTE$, let $X_{1:n}<\dots
<X_{n:n}$ denote the order statistics of $X_{1},\dots ,X_{n}$, we set%
\begin{equation}
\widetilde{\mathbb{C}}_{n}(t)=\dfrac{1}{1-t}\left[ \int_{t}^{1-k/n}\mathbf{Q}%
_{n}(s)ds+\int_{0}^{k/n}\mathbf{Q}_{n}^{w}(1-s)ds\right] ,  \label{eq 111}
\end{equation}%
where $\mathbf{Q}_{n}(s)$ is the empirical estimation of the
quantile function, which is equal to the i$^{th}$ order statistic $X_{i,n}$
for all $s\in \left( \left( i-1\right) /n,i/n\right) $ and for all $%
i=1,...,n $ and $\mathbf{Q}_{n}^{w}(1-s)$ is the estimation of the extreme quantile established by \cite{Weissman1978}, this estimator is given by%
\begin{equation}
\mathbf{Q}_{n}^{w}(1-s)=\left( \frac{k}{n}\right) ^{1/\widehat{\alpha }%
_{n}}X_{n-k,n}s^{-1/\widehat{\alpha }_{n}}.  \label{eq 6}
\end{equation}%
where $\widehat{\alpha }_{n}$  is  the Hill estimator \nocite{hill1975}\cite{hill1975}
\begin{equation}
\widehat{\alpha }_{n}=\left( {\frac{1}{k}}\sum_{i=1}^{k}\log
X_{n-i+1:n}-\log X_{n-k:n}\right) ^{-1}  \label{eq 7}
\end{equation}%
of the tail index $\alpha \in (1,2)$ and $k=k(n)$ is an intermediate integer
sequence satisfying the condition%
\begin{equation}
k\rightarrow \infty \text{ and }k/n\rightarrow 0\text{ as }n\rightarrow
\infty .  \label{eq 8}
\end{equation}%
After integration, we obtain
\begin{equation}
\widetilde{\mathbb{C}}_{n}(t)=\dfrac{1}{1-t}\left[ \int_{t}^{1-k/n}\mathbf{Q}%
_{n}(s)ds+\left( \frac{k}{n}\right) \frac{\widehat{\alpha }_{n}X_{n-k},_{n}}{%
\left( \widehat{\alpha }_{n}-1\right) }\right] .  \label{eq 9}
\end{equation}%
the asymptotic normality of this estimator is found in \cite{necir2010} under some conditions.

We note that, the Hill's estimator $\widehat{\alpha }_{n}$ plays a pivotal role in statistical
inference on distribution tails. This estimator has been extensively studied,
improved and even generalized to any real parameter $\alpha $. Weak
consistency of $\widehat{\alpha }_{n}$ was established by \cite{mason1982}
assuming only that the underlying cdf $F$ satisfies condition (\ref{eq1-1}). The
asymptotic normality of $\widehat{\alpha }_{n}$ has been established (see \cite{deHaanPeng1998}) under the following stricter condition that
characterizes Hall's model (see \cite{hall1982} and \cite{hall1985}),
where there exist $c>0$, $d\neq0$ and $\beta>\alpha>0$ such that%
\begin{equation}
1-F(x)=cx^{-\alpha }+dx^{-\beta }+o\left( x^{-\beta }\right) ,\text{ as }%
x\rightarrow \infty ,  \label{eq 10}
\end{equation}%
Note that (\ref{eq 10}), which is a special case of a more general
second-order regular variation condition (see \cite{deHaan1996}), is equivalent to%
\begin{equation}
\mathbf{Q}(1-s)=c^{1/\alpha }s^{-1/\alpha }(1+\alpha ^{-1}c^{-\beta /\alpha
}ds^{\beta /\alpha -1}+o(1)),\text{ as }s\downarrow 0,  \label{eq 11}
\end{equation}%
The constants $\alpha $ and $\beta $ are called, respectively, first-order
(tail index, shape parameter) and second-order parameters of cdf $F$.
In this paper, we use the bias-reduced estimator of the high quantile $%
\mathbf{Q}(1-s)$, proposed by \cite{li2010} who exploited the censored
maximum likelihood (CML) based estimators $\widehat{\alpha }$ and $\widehat{%
\beta }$ of the couple of regular variation parameters $(\alpha ,\beta )$
introduced by \cite{peng2004}. The CML estimators $\widehat{\alpha }$, $%
\widehat{\beta }$ are defined as the solution of the two equations (under
the constraint $\beta >\alpha ):$%
\begin{equation}
\frac{1}{k}\sum_{i=1}^{k}\frac{1}{G_{i}(\alpha ,\beta )}=1\text{ and }\frac{1%
}{k}\sum_{i=1}^{k}\frac{1}{G_{i}(\alpha ,\beta )}\log \frac{X_{n-i+1,n}}{%
X_{n-k,n}}=\beta ^{-1},  \label{eq 12}
\end{equation}%
where%
\begin{equation}
G_{i}(\alpha ,\beta )=\frac{\alpha }{\beta }\left( 1+\frac{\alpha \beta }{%
\alpha -\beta }H\left( \alpha \right) \right) \left( \frac{X_{n-i+1,n}}{%
X_{n-k,n}}\right) ^{\beta -\alpha }-\frac{\alpha \beta }{\alpha -\beta }%
H\left( \alpha \right) ,  \label{eq 13}
\end{equation}%
and%
\begin{equation}
H\left( \alpha \right) =\frac{1}{\alpha }-\frac{1}{k}\sum_{i=1}^{k}\log 
\frac{X_{n-i+1,n}}{X_{n-k,n}}.  \label{eq 14}
\end{equation}%
Li et al. \cite{li2010} obtained their bias-reduced estimators $\mathbf{Q}%
_{n}^{L}\left( 1-s\right) $, of the high quantiles $\mathbf{Q}(1-s)$, by
substituting $\left( \widehat{\alpha },\widehat{\beta }\right) $ to $(\alpha
,\beta )$ in (\ref{eq 11}). That is%
\begin{equation}
\widehat{\mathbf{Q}}_{n}^{L}(1-s)=\widehat{c}^{1/\widehat{\alpha }}s^{-1/%
\widehat{\alpha }}(1+\widehat{\alpha }^{-1}c^{-\widehat{\beta }/\widehat{%
\alpha }}\widehat{d}s^{\widehat{\beta }/\widehat{\alpha }-1}+o(1)),\text{ as 
}s\downarrow 0,  \label{eq 15}
\end{equation}%
where%
\begin{equation}
\left\{ 
\begin{array}{c}
\widehat{c}=\frac{\widehat{\alpha }\widehat{\beta }}{\widehat{\alpha }-%
\widehat{\beta }}\frac{k}{n}X_{n-k,n}^{\widehat{\alpha }}\left( \frac{1}{%
\widehat{\beta }}-\frac{1}{k}\sum_{i=1}^{k}\log \frac{X_{n-i+1,n}}{X_{n-k,n}}%
\right) \\ 
\widehat{d}=\frac{\widehat{\alpha }\widehat{\beta }}{\widehat{\alpha }-%
\widehat{\beta }}\frac{k}{n}X_{n-k,n}^{\widehat{\beta }}\left( \frac{1}{%
\widehat{\alpha }}-\frac{1}{k}\sum_{i=1}^{k}\log \frac{X_{n-i+1,n}}{X_{n-k,n}%
}\right)%
\end{array}%
\right. .  \label{eq 16}
\end{equation}%
The consistency and asymptotic normality of $\widehat{\mathbf{Q}}_{n}^{L}(1-s)$
are established by the same authors. By using $\mathbf{Q}_{n}$ in formula (%
\ref{eq 111}), we get, after integration, the new estimator of the $CTE$ as
follows:%
\begin{equation}
\overline{\mathbb{C}}_{n}(t)=\dfrac{1}{1-t}\left[ \int_{t}^{1-k/n}\mathbf{Q}%
_{n}(s)ds+\left( \frac{k}{n}\right) \left( \frac{n\widehat{c}}{k}\right) ^{1/%
\widehat{\alpha }}\left( \frac{\widehat{\alpha }}{\widehat{\alpha }-1}+\frac{%
\widehat{d} \widehat{c}^{-\widehat{\beta }/\widehat{\alpha }}\left( k/n\right)
^{\widehat{\beta }/\widehat{\alpha }-1}}{\widehat{\beta }-1}\right) \right] ,
\label{eq 17}
\end{equation}%
provided that $\widehat{\beta }>\widehat{\alpha }>1$ so that $\overline{%
\mathbb{C}}_{n}$ is finite.\medskip

The rest of the paper is organized as follows. In Section \ref{sect-2} we study the asymptotic normality of
the new $CTE$ estimator. In section \ref{simul}, we make a simulation study of the new estimator and we show the performance of our estimator and compare it with the old one. We finich our paper by a conclusion given in section \ref{conclusion}. The proof of the main result, which is Theorem \ref%
{th-new} in Section \ref{sect-2}, is postponed to Section \ref{sect-4}.

\section{Main result and its practical implementation}

\label{sect-2}In the field of the extreme values theory, a function denoted by $\mathbf{U}$
and (sometimes) called tail quantile function, is defined by 
\begin{equation}
\mathbf{U}(t)=F^{\leftarrow }(1-1/t)=\mathbf{Q}\left( 1-1/t\right) ,
\label{eq 2}
\end{equation}%
where $F^{\leftarrow }$ represented the generalized inverse  of the df $F.$
Then, we say that $F$ is heavy-tailed iff $\mathbf{U}\in
RV_{1/\alpha }$ (de Haan, 1970 \cite{haan1970}), i.e.
\begin{equation}
\underset{t\rightarrow \infty }{\lim }\frac{\mathbf{U}(tx)}{\mathbf{U}(t)}%
=x^{1/\alpha },\text{ for any }x>0.  \label{eq 3}
\end{equation}%
In terms of this function, Hall's conditions (\ref{eq 10}) and (\ref{eq 11})
are equivalent to
\begin{equation}
\mathbf{U}(t)=c^{1/\alpha }t^{1/\alpha }(1+\alpha ^{-1}c^{-\beta /\alpha
}dt^{1-\beta /\alpha }+o(1)),\text{ as }t\rightarrow \infty .  \label{eq 11'}
\end{equation}%
This implies that, there is a function $A_{1}(t)$, which tends to zero as $%
t\rightarrow \infty $ (because $\beta >\alpha $), determines the rate of
convergence of $\log \left( \mathbf{U}(tx)/\mathbf{U}(t)\right) $ to its
limit $\alpha ^{-1}\log x$, such that:%
\begin{equation}
\underset{t\rightarrow \infty }{\lim }\frac{\log \left( \mathbf{U}(tx)/%
\mathbf{U}(t)\right) -\alpha ^{-1}\log x}{A_{1}(t)}=\frac{x^{1-\beta /\alpha
}-1}{1-\beta /\alpha };\text{ for any }x>0,  \label{eq 20}
\end{equation}%
where%
\begin{equation*}
A_{1}(t)=d\alpha ^{-1}\left( 1-\beta /\alpha \right) c^{-\beta /\alpha
}t^{1-\beta /\alpha }.
\end{equation*}%
Relation (\ref{eq 20}) is known as the second-order condition of regular
variation (see, e.g., \cite{deHaan2006}). Unfortunately,
the second-order regular variation is not sufficient to find asymptotic
distributions for the estimators defined by the systems (\ref{eq 12}) and (%
\ref{eq 16}). We strengthen it with a condition, called third-order
condition of regular variation and given by (\ref{eq 21}), that specifies
the rate of (\ref{eq 20}) (see, e.g., \cite{deHaan1996} or
\cite{alves2007}).%
\begin{equation}
\underset{t\rightarrow \infty }{\lim }\frac{\frac{\log \left( \mathbf{U}(tx)/%
\mathbf{U}(t)\right) -\alpha ^{-1}\log x}{A_{1}(t)}-\frac{x^{1-\beta /\alpha
}-1}{1-\beta /\alpha }}{A_{2}(t)}=D\left( \alpha ,\beta ,\rho \right) \frac{%
x^{1-\beta /\alpha }-1}{1-\beta /\alpha };\text{ for any }x>0.  \label{eq 21}
\end{equation}%
where $A_{2}(t)\rightarrow 0$ as $t\rightarrow \infty $, with constant sign
near infinity and%
\begin{equation*}
D\left( \alpha ,\beta ,\rho \right) =\frac{1}{\rho }\left( \frac{x^{1-\beta
/\alpha +\rho }-1}{1-\beta /\alpha +\rho }-\frac{x^{1-\beta /\alpha }-1}{%
1-\beta /\alpha }\right) .
\end{equation*}%
with $\rho $ being a positive constant called third-order parameter. Peng
and Qi, \cite{peng2004} established the asymptotic normality of $\widehat{\alpha }$, $%
\widehat{\beta }$ and $\widehat{c}$ under the following extract conditions
on the sample fraction $k$, as $n\rightarrow \infty $:
\begin{equation}
(i):k^{1/2}\left\vert A_{1}(n/k)\right\vert \rightarrow \infty ,\text{ }%
(ii):k^{1/2}A_{1}^{2}(n/k)\rightarrow 0\text{ and }%
(iii):k^{1/2}A_{1}(n/k)A_{2}(n/k)\rightarrow 0.  \label{eq 22}
\end{equation}%
As for $\widehat{d}$, it is asymptotically normal under the assumption $%
k^{1/2}\left\vert A_{1}(n/k)\right\vert /\log (n/k)\rightarrow \infty $
added to $(ii)$ and $(iii)$. Note that, from a theoretical point of view,
assumptions (\ref{eq 8}) and (\ref{eq 22}) are realistic, as the following
example shows, indeed, let us choose 
\begin{equation}
k=[n^{1-\varepsilon }],0<\varepsilon <1,  \label{eq 23}
\end{equation}%
then it easy to verify that these assumptions hold for any $1/5<\varepsilon
<1/3$. The notation [\textperiodcentered ] stands for the integer part of
real numbers.\medskip

Our main result, namely the asymptotic normality of new estimator of the
$CTE$ is formulated in the following theorem.

\begin{theorem}\label{th-new}
Assume that the cdf $F$ satisfies condition (\ref{eq 21})
with $\alpha \in (1,2)$ and $\beta /\alpha =\lambda >1$. Then for any
sequence of integers $k=k_{n}$ satisfying the conditions (\ref{eq 8}) and (\ref{eq 22}). Then
\begin{equation}
\frac{\sqrt{n}\big(\overline{\mathbb{C}}_{n}(t)-\mathbb{C}(t)\big)(1-t)}{%
(k/n)^{1/2}\left( nc/k\right) ^{1/\alpha }}\rightarrow _{d}\mathcal{N}\left(
0,\sigma ^{2}\left( \alpha ,\beta \right) \right)  \label{statement-new}
\end{equation}%
for any fixed $t\in (0,1)$, where the asymptotic variance $\sigma ^{2}\left(
\alpha ,\beta \right) $ is given by the formula 
\begin{equation*}
\sigma ^{2}\left( \alpha ,\beta \right) ={\frac{\alpha ^{2}\beta ^{4}}{%
(\alpha -1)^{4}(\alpha -\beta )^{4}}+}\frac{2}{2-\alpha }+\frac{2\alpha \beta
^{2}}{(\alpha -1)^{2}(\alpha -\beta )^{2}}.
\end{equation*}
\end{theorem}
\section{Simulation study}
\label{simul}
In this section, the biased estimator $\tilde{\mathbb{C}}_{n}(t)$ and the reduced-bias one  $
\overline{\mathbb{C}}_{n}(t)$ are compared using simulation study. For this reason, $1000$ samples of size $n \in \left\lbrace 250, 500, 1000, 2000\right\rbrace $ are simulated
from a two heavy-tailed distributions :
\begin{description}
\item \textbf{Fr\'{e}chet model :} defined as :  $ F(x)=exp(-x^{-\alpha}) $  we take   $\alpha =1,5$ and $\alpha =1,75,$ respectively,  and according  to the Hall's model, we find $\beta =2\alpha , c=1,$ and $d=-1/2$. 
\item \textbf{Burr model :} defined as : $ F(x)=1-(1+x^\tau)^{-\lambda} $  with $\lambda =1,5$ and $\lambda =1,75,$ both with $\tau=1$ and in this case we have $ c=1,\alpha=\tau\lambda,d=-\lambda $ and $ \beta=\lambda\tau+\tau$ ( For more flexibility, we can take $ \tau=\alpha/\lambda $ in the expression of Burr distribution).
\end{description}
In all cases, we assume that $\alpha ,\beta ,c$ and $d $ are unknowns, then by the resolving the equations (\ref{eq 12}) and (\ref{eq 13}) we calculate the values of $\widehat{\alpha }$ and $\widehat{\beta }$, and resolving the
system (\ref{eq 16}) to calculate the estimators $\widehat{c}$ and $%
\widehat{d}.$ We fix two values of $t,$ for example $t=0.90$ and $t=0.95$ to calculate the differents values of the two estimators, the new estimator $\overline{\mathbb{C}}_{n}$ and the old estimator $\widetilde{\mathbb{C}}_{n}$.  We note that, the mean, the bias  and root mean squared error (RMSE) of these estimators are estimated over the $1000$ replications.
We illustrate and compare  the bias  and the root mean square error (RMSE) of  $\tilde{\mathbb{C}}_{n}(t)$ and  $\overline{\mathbb{C}}_{n}(t)$. Finally,  we  summarize the results of simulations in the tables (\ref{table1}, \ref{table2}, \ref{table3}, and \ref{table4}). We remark, on one hand, that the bias and the RMSE of the new estimator increase when the sample size increase, and on other hands, the bias and RMSE of our estimator is smaller than the old estimator.\\
In order to show the influence of the choice of the integer value $k$ to  the performance of our estimator $\overline{\mathbb{C}}_{n}(t)$  and the old estimator $\tilde{\mathbb{C}}_{n}(t)$, we generate $  1000$ samples of size $1000$ of the parent model ( Fr\'{e}chet and Burr) with two values of index $\alpha=1.5$ and $1.75$ and we picture the two estimators when $k$ varied from $50$ to $850$. The results are displayed in figures \ref{figure1} and \ref{figure2}. We observe that the performance of our estimator is clearly better than the old estimator.
\begin{center}
\begin{table}[H]
\caption{Comparison between the new estimator $\overline{\mathbb{C}}_{n}$ and old estimator $\widetilde{\mathbb{C}}_{n}$  of $CTE$ in terms of bias and RMSE  respecting to the variation of sample size $ n $ based on Fr\'{e}chet distribution.}
\begin{equation*}
\begin{tabular}{|c|c|c|c|c|c|c|c|c|}
\hline
\multicolumn{9}{|c|}{Fr\'{e}chet distribution, $\alpha =1.5$} \\ \hline\hline
$t$ & $n$ & $CTE$ & $\mathbb{\bar{C}}_{n}\left( t\right) $ & bias & RMSE & $%
\mathbb{\tilde{C}}_{n}\left( t\right) $ & bias & RMSE \\ \hline
$0.9$ & $250$ & $13.793$ & $14.457$ & $-0.664$ & $0.66415$ & $11.048$ & $%
2.745$ & $2.7454$ \\ \cline{2-2}\cline{4-9}
& $500$ &  & $14.409$ & $-0.616$ & $0.61578$ & $11.418$ & $2.375$ & $2.3755$
\\ \cline{2-2}\cline{4-9}
& $1000$ &  & $14.287$ & $-0.494$ & $0.49345$ & $11.837$ & $1.956$ & $1.9559$
\\ \cline{2-2}\cline{4-9}
& $2000$ &  & $14.14$ & $-0.347$ & $0.34733$ & $11.852$ & $1.941$ & $1.941$
\\ \hline\hline
$0.95$ & $250$ & $21.984$ & $22.78$ & $-0.796$ & $0.79685$ & $16.115$ & $%
-5.869$ & $5.8683$ \\ \cline{2-2}\cline{4-9}
& $500$ &  & $22.703$ & $-0.719$ & $0.71965$ & $16.473$ & $-5.511$ & $5.5105$
\\ \cline{2-2}\cline{4-9}
& $1000$ &  & $22.674$ & $-0.69$ & $0.69043$ & $16.944$ & $-5.04$ & $5.0392$
\\ \cline{2-2}\cline{4-9}
& $2000$ &  & $22.67$ & $-0.686$ & $0.68661$ & $17.247$ & $-4.737$ & $4.7364$
\\ \hline
\end{tabular}
\end{equation*}
\label{table1}
\end{table}
\begin{table}[H]
\caption{Comparison between the new estimator $\overline{\mathbb{C}}_{n}$ and old estimator $\widetilde{\mathbb{C}}_{n}$  of $CTE$ in terms of bias and RMSE  respecting to the variation of sample size $ n $ based on Fr\'{e}chet distribution.}
\begin{equation*}
\begin{tabular}{|c|c|c|c|c|c|c|c|c|}
\hline
\multicolumn{9}{|c|}{Fr\'{e}chet distribution, $\alpha =1.75$} \\ \hline\hline
$t$ & $n$ & $CTE$ & $\mathbb{\bar{C}}_{n}\left( t\right) $ & bias & RMSE & $%
\mathbb{\tilde{C}}_{n}\left( t\right) $ & bias & RMSE \\ \hline
$0.9$ & $250$ & $8.6207$ & $8.6977$ & $-0.077$ & $0.076993$ & $8.2478$ & $%
0.3729$ & $0.37295$ \\ \cline{2-2}\cline{4-9}
& $500$ &  & $8.6771$ & $-0.0564$ & $0.05637$ & $8.3136$ & $0.3071$ & $%
0.30711$ \\ \cline{2-2}\cline{4-9}
& $1000$ &  & $8.6586$ & $-0.0379$ & $0.037901$ & $8.5085$ & $0.1122$ & $%
0.11218$ \\ \cline{2-2}\cline{4-9}
& $2000$ &  & $8.6504$ & $-0.0297$ & $0.029695$ & $8.5374$ & $0.0833$ & $%
0.083322$ \\ \hline\hline
$0.95$ & $250$ & $12.866$ & $12.963$ & $-0.097$ & $0.096298$ & $12.185$ & $%
0.681$ & $0.68126$ \\ \cline{2-2}\cline{4-9}
& $500$ &  & $12.956$ & $-0.09$ & $0.089431$ & $12.218$ & $0.648$ & $0.64858$
\\ \cline{2-2}\cline{4-9}
& $1000$ &  & $12.937$ & $-0.071$ & $0.070963$ & $12.624$ & $0.242$ & $%
0.24248$ \\ \cline{2-2}\cline{4-9}
& $2000$ &  & $12.917$ & $-0.051$ & $0.050889$ & $12.766$ & $0.100$ & $%
0.10069$ \\ \hline
\end{tabular}%
\end{equation*}
\label{table2}
\end{table}
\begin{table}[H]
\caption{Comparison between the new estimator $\overline{\mathbb{C}}_{n}$ and old estimator $\widetilde{\mathbb{C}}_{n}$  of $CTE$ in terms of bias and RMSE  respecting to the variation of sample size $ n $ based on Burr distribution.}
\begin{equation*}
\begin{tabular}{|c|c|c|c|c|c|c|c|c|}
\hline
\multicolumn{9}{|c|}{Burr distribution, $\alpha =1.5$} \\ 
\hline\hline
$t$ & $n$ & $CTE$ & $\mathbb{\bar{C}}_{n}\left( t\right) $ & bias & RMSE & $%
\mathbb{\tilde{C}}_{n}\left( t\right) $ & bias & RMSE \\ \hline
$0.9$ & $250$ & $13.676$ & $14.365$ & $-0.689$ & $0.68807$ & $14.909$ & $%
-1.233$ & $1.2321$ \\ \cline{2-2}\cline{4-9}
& $500$ &  & $14.226$ & $-0.55$ & $0.54908$ & $14.658$ & $-0.982$ & $0.98132$
\\ \cline{2-2}\cline{4-9}
& $1000$ &  & $14.124$ & $-0.448$ & $0.44795$ & $14.506$ & $-0.83$ & $0.82922
$ \\ \cline{2-2}\cline{4-9}
& $2000$ &  & $14.026$ & $-0.35$ & $0.34997$ & $14.491$ & $-0.815$ & $0.81428
$ \\ \hline\hline
$0.95$ & $250$ & $21.891$ & $23.239$ & $-1.348$ & $1.3474$ & $24.05$ & $%
-2.159$ & $2.1587$ \\ \cline{2-2}\cline{4-9}
& $500$ &  & $22.996$ & $-1.105$ & $1.1044$ & $23.874$ & $-1.983$ & $1.9827$
\\ \cline{2-2}\cline{4-9}
& $1000$ &  & $22.753$ & $-0.862$ & $0.86155$ & $23.709$ & $-1.818$ & $1.8173
$ \\ \cline{2-2}\cline{4-9}
& $2000$ &  & $22.602$ & $-0.711$ & $0.71095$ & $23.64$ & $-1.749$ & $1.7489$
\\ \hline
\end{tabular}
\end{equation*}
\label{table3}
\end{table}
\begin{table}[H]
\caption{Comparison between the new estimator $\overline{\mathbb{C}}_{n}$ and old estimator $\widetilde{\mathbb{C}}_{n}$  of $CTE$ in terms of bias and RMSE  respecting to the variation of sample size $ n $ based on Burr distribution.}
\begin{equation*}
\begin{tabular}{|c|c|c|c|c|c|c|c|c|}
\hline
\multicolumn{9}{|c|}{Burr distribution, $\alpha =1.75$} \\ 
\hline\hline
$t$ & $n$ & $CTE$ & $\mathbb{\bar{C}}_{n}\left( t\right) $ & bias & RMSE & $%
\mathbb{\tilde{C}}_{n}\left( t\right) $ & bias & RMSE \\ \hline
$0.9$ & $250$ & $8.5455$ & $9.0332$ & $-0.4877$ & $0.48776$ & $9.3583$ & $%
-0.8128$ & $0.81281$ \\ \cline{2-2}\cline{4-9}
& $500$ &  & $8.9188$ & $-0.3733$ & $0.37334$ & $9.0502$ & $-0.5047$ & $%
0.5047$ \\ \cline{2-2}\cline{4-9}
& $1000$ &  & $8.8353$ & $-0.2898$ & $0.28979$ & $8.9942$ & $-0.4487$ & $%
0.44873$ \\ \cline{2-2}\cline{4-9}
& $2000$ &  & $8.7708$ & $-0.2253$ & $0.22533$ & $8.9115$ & $-0.366$ & $%
0.36607$ \\ \hline\hline
$0.95$ & $250$ & $12.811$ & $13.697$ & $-0.886$ & $0.88641$ & $14.386$ & $%
-1.575$ & $1.5751$ \\ \cline{2-2}\cline{4-9}
& $500$ &  & $13.556$ & $-0.745$ & $0.74521$ & $14.082$ & $-1.271$ & $1.2709$
\\ \cline{2-2}\cline{4-9}
& $1000$ &  & $13.381$ & $-0.57$ & $0.56966$ & $13.788$ & $-0.977$ & $0.97676
$ \\ \cline{2-2}\cline{4-9}
& $2000$ &  & $13.246$ & $-0.435$ & $0.4355$ & $13.596$ & $-0.785$ & $0.78538
$ \\ \hline
\end{tabular}%
\end{equation*}
\label{table4}
\end{table}
\begin{figure}[H]
	\centering
	\includegraphics[width=1.0\linewidth]{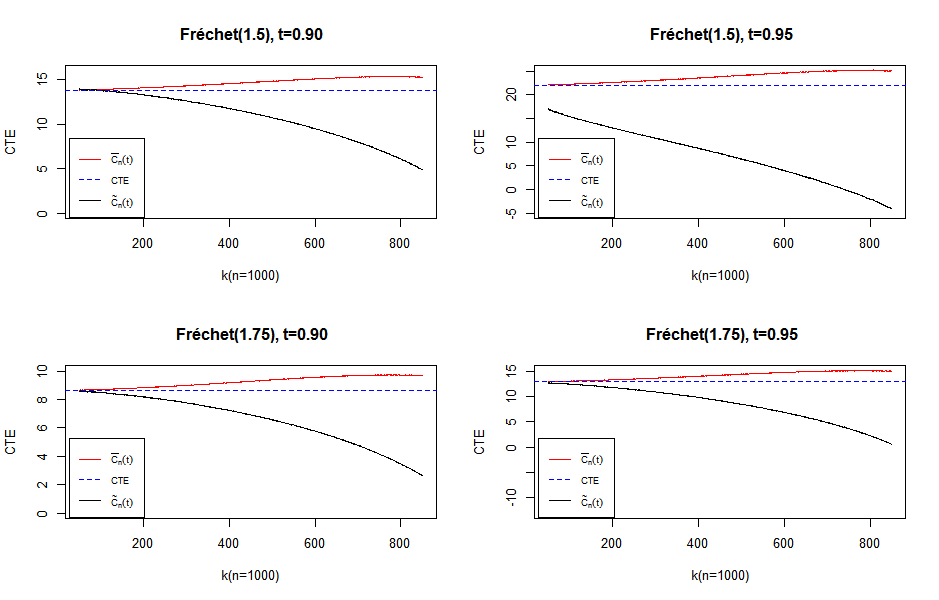}
	\caption{Behavior of the new estimator $\overline{\mathbb{C}}_{n}$ and old estimator $\widetilde{\mathbb{C}}_{n}$  of $CTE$ in respect to the variation of $ k $ for Fr\'{e}chet distribution.}
	\label{figure1}
\end{figure}

\begin{figure}[H]
	\centering
	\includegraphics[width=0.9\linewidth]{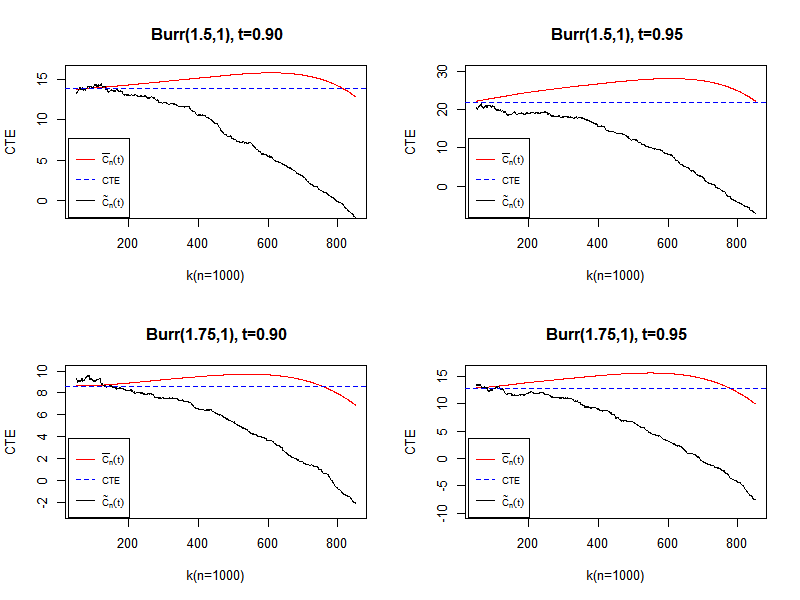}
	\caption{Behavior of the new estimator $\overline{\mathbb{C}}_{n}$ and old estimator $\widetilde{\mathbb{C}}_{n}$  of $CTE$ in respect to the variation of $ k $ for Burr distribution.}
	\label{figure2}
\end{figure}
\end{center}
The analysis and simulation as shown in tables \ref{table1},  \ref{table2}, \ref{table3}, \ref{table4}, and figures \ref{figure1}, \ref{figure2} indicate that the new $CTE$ estimator is more performant than the estimator contructed by Necir et al (2010) \cite{necir2010}.
\section{Conclusion}
\label{conclusion}
From the research that has undertaken, it is possible to conclude that the estimation of the $CTE$ by the use of extreme quantile established by Li (2010) \cite{li2010} gave us a more efficient estimate compared to that introduced by Necir et al (2010) \cite{necir2010}. the results of the simulations show significant improvement either in term of bias or in term of RMSE.\\ The next stage of our research will be the construction of a confidence interval of our $CTE$ estimator with the covering probability. 
\section{Proof of Theorem \protect\ref{th-new}}
\label{sect-4}
\begin{flushleft}
We start the proof of Theorem \ref{th-new} with the decomposition%
\begin{equation*}
\left( 1-t\right) \mathbb{C}(t)=C_{1.n}(s)+C_{2,n}(s),
\end{equation*}%
where%
\begin{equation*}
C_{1.n}(s)=\int_{t}^{1-k/n}\mathbf{Q}\left( s\right) ds\text{ and }%
C_{2.n}(s)=\int_{0}^{k/n}\mathbf{Q}\left( 1-s\right) ds.
\end{equation*}%
Also, we have 
\begin{equation*}
\left( 1-t\right) \overline{\mathbb{C}}(t)=\overline{C}_{1.n}(s)+\overline{C}%
_{2,n}(s),
\end{equation*}%
where%
\begin{equation*}
\overline{C}_{1.n}(s)=\int_{t}^{1-k/n}\mathbf{Q}_{n}(s)ds,
\end{equation*}%
and%
\begin{equation*}
\overline{C}_{2.n}(s)=\left( \frac{k}{n}\right) \left( \frac{n\widehat{c}}{k}%
\right) ^{1/\widehat{\alpha }}\left( \frac{\widehat{\alpha }}{\widehat{%
\alpha }-1}+\frac{\widehat{d}\widehat{c}^{-\widehat{\beta }/\widehat{\alpha }%
}\left( k/n\right) ^{\widehat{\beta }/\widehat{\alpha }-1}}{\widehat{\beta }%
-1}\right) .
\end{equation*}%
First, it is easy to verify that, as $n\rightarrow \infty $%
\begin{equation*}
C_{2.n}(s)=\left( 1+o_{p}\left( 1\right) \right) \left( \frac{k}{n}\right)
\left( \frac{nc}{k}\right) ^{1/\alpha }\left( \frac{\alpha }{\alpha -1}%
\right) ,
\end{equation*}%
and, under the condition (\ref{eq 7}), we have%
\begin{equation*}
\overline{C}_{2.n}(s)=\left( 1+o_{p}\left( 1\right) \right) \left( \frac{k}{n%
}\right) \left( \frac{n\widehat{c}}{k}\right) ^{1/\widehat{\alpha }}\left( 
\frac{\widehat{\alpha }}{\widehat{\alpha }-1}\right) .
\end{equation*}%
It follows that 
\begin{equation}
\overline{C}_{2.n}(s)-C_{2.n}(s)=\left( 1+o_{p}\left( 1\right) \right)
\left( \frac{k}{n}\right) \left[ \frac{\widehat{\alpha }}{\widehat{\alpha }-1%
}\left( \frac{n\widehat{c}}{k}\right) ^{1/\widehat{\alpha }}-\frac{\alpha }{%
\alpha -1}\left( \frac{nc}{k}\right) ^{1/\alpha }\right] ,  \label{stat-0}
\end{equation}%
let us write $\overline{C}_{2.n}(s)-C_{2.n}(s)=A_{n,1}(s)+A_{n,2}(s),$ where%
\begin{gather*}
A_{n,1}(t)=\left( 1+o_{p}\left( 1\right) \right) \sqrt{k}\left( \frac{%
\widehat{\alpha }}{\widehat{\alpha }-1}-\frac{\alpha }{\alpha -1}\right) ,
A_{n,2}=\left( 1+o_{p}\left( 1\right) \right) \sqrt{k}\left( \frac{\left( n%
\widehat{c}/k\right) ^{1/\widehat{\alpha }}}{\left( nc/k\right) ^{1/\alpha }}%
-1\right) .
\end{gather*}%
We begin by showing that $A_{n,2}\rightarrow 0$, as $n\rightarrow \infty $.
First observe that $A_{n,2}$ may be rewritten into%
\begin{equation*}
A_{n,2}=\left( 1+o_{p}\left( 1\right) \right) \sqrt{k}\left( \left(
nc/k\right) ^{1/\widehat{\alpha }-1/\alpha }\left( \widehat{c}/c\right) ^{1/%
\widehat{\alpha }-1/\alpha }-1\right)
\end{equation*}%
Assumptions $(i)$ and $(ii)$, imply that $\sqrt{k}/\log (n/k)\rightarrow
\infty $. Also, from Theorem 1 of \cite{peng2004}, the asymptotic
normality of $\widehat{\alpha }$ gives $\widehat{\alpha }-\alpha =O\left( 1/%
\sqrt{k}\right) .$ Therefore%
\begin{equation*}
\left( 1/\widehat{\alpha }-1/\alpha \right) \log \left( nc/k\right) \overset{%
P}{\rightarrow }0,
\end{equation*}%
this implies that%
\begin{equation*}
\left( nc/k\right) ^{1/\widehat{\alpha }-1/\alpha }\overset{P}{\rightarrow }1%
\text{ as }n\rightarrow \infty \text{.}
\end{equation*}%
On the other hand, from equation (4.7) in \cite{li2010}, we have%
\begin{equation*}
\left( \widehat{c}/c\right) -1=\alpha ^{-1}\left( 1+o_{p}\left( 1\right)
\right) \left( \widehat{\alpha }-\alpha \right) \log (n/k)+o_{p}\left( \frac{%
1}{\sqrt{k}}\log \frac{n}{k}\right) .
\end{equation*}%
Since $\widehat{c}$ is a consistent estimator of $c$, then Taylor's
expansion gives%
\begin{equation*}
\left( \widehat{c}/c\right) ^{1/\widehat{\alpha }-1/\alpha }-1=\alpha
^{-1}\left( \widehat{\alpha }-\alpha \right) \log (n/k)\left( \widehat{c}%
/c-1\right) ,\text{ as }n\rightarrow \infty .
\end{equation*}%
It suffices now to show that $\sqrt{k}\left( \left( \widehat{c}/c\right) ^{1/%
\widehat{\alpha }-1/\alpha }-1\right) $ converges to $0$ in probability.
Indeed, again by using the fact that $\widehat{\alpha }-\alpha =O_{p}\left(
1/\sqrt{k}\right) ,$ yield%
\begin{equation*}
\sqrt{k}\left( \left( \widehat{c}/c\right) ^{1/\widehat{\alpha }-1/\alpha
}-1\right) =O_{p}\left( 1\right) \left( \frac{1}{\sqrt{k}}\log \frac{n}{k}%
+o_{p}\left( \frac{\log (n/k)}{\sqrt{k}}\right) \right) ,
\end{equation*}%
which tends in probability to $0$, because we already have $\sqrt{k}/\log
(n/k)\rightarrow \infty .$ Now, we consider the term $A_{n,1}$. Since $%
\widehat{\alpha }$ is a consistent estimator of $\alpha $, then it is easy
to show that%
\begin{equation*}
A_{n,1}=-\left( 1+o_{p}\left( 1\right) \right) \frac{\sqrt{k}}{\left( \alpha
-1\right) ^{2}}\left( \widehat{\alpha }-\alpha \right) .
\end{equation*}%
From theorem 2.1 in \cite{brahimi2013}, we infer that%
\begin{equation*}
\frac{\sqrt{n}\left( 1-t\right) }{(k/n)^{1/2}\left( nc/k\right) ^{1/\alpha }}%
\left( \overline{C}_{2,n}(s)-C_{2,n}(s)\right) =-\frac{\alpha }{\left(
\alpha -1\right) ^{2}}\left[ \eta _{1}W_{1}+\eta _{2}W_{2}+\eta _{3}W_{3}%
\right] +o_{\mathbf{P}}(1),
\end{equation*}
where
\begin{eqnarray*}
W_{1} &=&\sqrt{n/k}B_{n}\left( 1-k/n\right) -\sqrt{n/k}%
\int_{0}^{1}s^{-1}B_{n}(1-ks/n)ds, \\
W_{2} &=&\left( \lambda ^{-1}-1\right) \sqrt{n/k}B_{n}\left( 1-k/n\right)
+\left( \lambda -1\right) \sqrt{n/k}\int_{0}^{1}s^{\lambda
-2}B_{n}(1-ks/n)ds, \\
W_{3} &=&\left( 1-\lambda \right) \sqrt{n/k}\int_{0}^{1}s^{\lambda -2}(\log
s)B_{n}(1-ks/n)ds \\
&&+\lambda ^{-2}\sqrt{n/k}B_{n}\left( 1-k/n\right) -\sqrt{n/k}%
\int_{0}^{1}s^{\lambda -2}B_{n}(1-ks/n)ds,
\end{eqnarray*}%
and%
\begin{equation*}
\eta _{1}=\frac{\lambda ^{4}}{\left( \lambda -1\right) ^{4}},\eta _{2}=\frac{%
\lambda ^{2}\left( 2\lambda -1\right) \left( 3\lambda -1\right) }{\left(
\lambda -1\right) ^{5}},\eta _{3}=\frac{\lambda ^{3}\left( 2\lambda
-1\right) ^{2}}{\left( \lambda -1\right) ^{4}}.
\end{equation*}%
From the Proof of statement (4.2) in \cite{necir2010}, we shall show
below that there are Brownian bridges $B_{n}$ such that 
\begin{equation}
\frac{\sqrt{n}\left( 1-t\right) }{(k/n)^{1/2}\mathbf{Q}(1-k/n)}\left( \overline{C}%
_{1,n}(s)-C_{1,n}(s)\right) =-\frac{\int_{k/n}^{1}B_{n}(1-s)d\mathbf{Q}(1-s)}{%
(k/n)^{1/2}\mathbf{Q}(1-k/n)}+o_{\mathbf{P}}(1)  \label{stat-1}
\end{equation}%
On the other hand, from (\ref{eq 11}), we have $\mathbf{Q}(1-k/n)\sim \left(
nc/k\right) ^{1/\alpha }$ , as $n\rightarrow \infty $, it follows that%
\begin{equation*}
\frac{\sqrt{n}\left( 1-t\right) }{(k/n)^{1/2}\left( nc/k\right) ^{1/\alpha }}%
\left( \overline{C}_{1,n}(s)-C_{1,n}(s)\right) =W_{4}+o_{\mathbf{P}}(1),
\end{equation*}%
where%
\begin{equation*}
W_{4}=-\frac{\int_{k/n}^{1}B_{n}(1-s)d\mathbf{Q}(1-s)}{(k/n)^{1/2}\mathbf{Q}(1-k/n)},
\end{equation*}
Finally, we have
\begin{equation*}
\frac{\sqrt{n}\left( 1-t\right) }{(k/n)^{1/2}\left( nc/k\right) ^{1/\alpha }}%
\left( \overline{\mathbb{C}}_{n}(t)-\mathbb{C}(t)\right) =-\frac{\alpha }{%
\left( \alpha -1\right) ^{2}}\left[ \eta _{1}W_{1}+\eta _{2}W_{2}+\eta
_{3}W_{3}\right] +W_{4}+o_{\mathbf{P}}(1).
\end{equation*}%
All done, after a standard calculation we obtain%
\begin{equation*}
\frac{\sqrt{n}\left( 1-t\right) }{(k/n)^{1/2}\left( nc/k\right) ^{1/\alpha }}%
\left( \overline{\mathbb{C}}_{n}(t)-\mathbb{C}(t)\right) \rightarrow
_{d}N\left( 0,\sigma ^{2}\left( \alpha ,\beta \right) \right) .
\end{equation*}
\end{flushleft}

\end{document}